\documentclass[10pt]{amsart}
\usepackage{latexsym}
\usepackage{amssymb}
\usepackage{amsmath}
\begin{document}

\title{
THE CALL OF MATHEMATICS
}
\author{S. S. Kutateladze}
\date{May 27, 2007}
\address[]{
Sobolev Institute of Mathematics\newline
\indent 4 Koptyug Avenue\newline
\indent Novosibirsk, 630090\newline
\indent RUSSIA
}
\begin{abstract}
A few remarks on  how mathematics quests for freedom.
\end{abstract}
\email{
sskut@member.ams.org
}
\maketitle


\medskip
Mathematics prevails in knowledge
as  the most ancient of sciences. However, in the beginning
was the  word. We must remember that the olden ``logos'' resides
beyond grammar. Today's mathematics became
the bastion of  logic, the savior of the order of mind and
the objectivity of  reasoning.

The intellectual field  resides beyond the grips of the law of
diminishing returns. The more we know, the huger become the frontiers
with the unbeknown, the oftener we meet the mysterious. The twentieth
century enriched our geometrical views with the concepts of space-time
and fractality.  Each instance of knowledge is an event, a point in the
Minkowski 4-space. The realm of our knowledge comprises a clearly bounded
set of these instances. The frontiers of science produce the boundary between the
known and the unknown which is undoubtedly fractal and we have no grounds
to assume it rectifiable or measurable. It is worth noting in parentheses that rather smooth are the routes to the
frontiers of science which are charted by teachers, professors, and all other kinds
of educationalists.
Pedagogics  dislikes saltations and sharp changes of the prevailing
paradigm. Possibly, these topological obstructions reflect some objective
difficulties in modernizing education.
The proofs are uncountable of the fractality of the
boundary between the known and the unbeknown. Among them we
see such negative trends as the unleashed growth of pseudoscience, mysticism,
and other forms of obscurantism which creep into all lacunas of the
unbeknown.   As revelations of fractality appear the most unexpected,
beautiful, and stunning interrelations between
seemingly distant areas and directions of
science. Mathematics serves as the principal catalyst of the unity of science.
There is evidence  galore  of the indispensability of
mathematics  in modernization and sustainable development.

We are granted the blissful world that has the indisputable property
of unique existence. The solitude of reality was perceived by our ancestors as
the ultimate proof of unicity.
Mathematics has never liberated itself from the tethers of
experimentation. The reason is not the simple fact that we still
complete proofs by declaring ``obvious.'' Alive and rather popular are
the views of mathematics as a toolkit for the natural sciences. These
stances may be expressed by the slogan ``mathematics is experimental
theoretical physics.'' Not less popular is the dual claim
``theoretical physics is experimental mathematics.'' This coupled
mottoes reflect the close affinity of
the trails of thought
in mathematics and the natural sciences.

It is worth observing that the dogmata of faith and the principles of
theology are also well reflected in the history of mathematical
theories. Variational calculus was invented in search of better
understanding of the principles of mechanics, resting on the religious
views of the universal beauty and harmony of the act of creation.

Mathematics is a rather specific area of intellectual
creativity which possess its own unparallel particularities.
Georg Cantor, the founder of set theory,
wrote in one of his classical papers in 1883 as follows:
``...das {\it Wesen\/} der {\it Mathematik\/}
liegt grerade in ihrer {\it Freiheit}.''
In other words, ``the essence of mathematics resides in its freedom.''
The freedom of modern mathematics does not reduce to
the absence of exogenous limitations of the objects and methods of research.
To a great extent, the freedom of mathematics is disclosed in
the new intellectual tools it provides for taming the universe,
liberating  a human being, and expanding the boundaries of his or her
independence.

The twentieth century marked an important twist in the content of
mathematics. Mathematical ideas imbued the humanitarian sphere and,
primarily, politics, sociology, and economics. Social events are
principally volatile and possess a high degree of uncertainty.
Economic processes utilize a wide range of the admissible ways of
production, organization, and management. The nature of nonunicity in
economics transpires: The genuine interests of human beings cannot
fail to be contradictory. The unique solution is an oxymoron in any
nontrivial problem of economics which refers to the distribution of
goods between a few agents. It is not by chance that the social
sciences and instances of humanitarian mentality invoke the numerous
hypotheses of the best organization of production and consumption, the
justest social structure, the codices of rational
behavior and moral conduct, etc.

The twentieth century became the age of freedom. Plurality and unicity were
confronted as collectivism and individualism. Many particular
phenomena of life and culture reflect their distinction. The
dissolution of monarchism and tyranny were accompanied by the rise of
parliamentarism and democracy.
In mathematics the quest for plurality led to the abandonment of the
overwhelming pressure of unicity and categoricity. The latter ideas were
practically absent, at least minor, in Ancient Greece and sprang to
life in the epoch of absolutism and Christianity.
Quantum mechanics and  Heisenberg's uncertainty incorporated plurality in physics.
The waves of modernism
in poetry and artistry should be also listed. Mankind had changed all
valleys of residence and dream.

The thesis of universal mathematization
enlightens many trends of today's thought.
Many  new synthetical areas of research
are the gains of  mathematics  which are decorated with
outstanding advances in economical cybernetics, theoretical programming,
mathematical linguistics, mathematical chemistry, and mathematical biology.
Mathematization of  the human sciences and the human dimension of
the natural sciences are familiar features of modernism.

Mathematics is a human science involving
the abstractions in which the human beings perceive  forms and relations.
Mathematics is impossible without the disciples,
professional mathematicians. Obviously, the essence of mathematics is disclosed
to us only as  expressed  in the contributions of scientists.
Therefore, it would be not a great exaggeration to paraphrase the
words  of Cantor and say that
{\it the essence of a mathematician resides and reveals itself in his or her
freedom}.

In science
we appraise and appreciate that which makes us wiser.
The notions of a good theory open up new possibilities
of solving particular problems. Rewarding is
the problem whose solution paves way to new fruitful concepts and
methods. Condescension is the mother of mediocrity. A fresh product of
a~mediocrity is called a banality. Time makes banal the most
splendid achievements, seminal theories, and challenging problems.
Indispensability is the most important quality of a good problem or
theory which refrains us from producing  banalities.

The greatest minds create indispensable scientific concepts and
ponder them over. They pose indispensable scientific problems and
contemplate over their solutions. The indispensable theories
and problems propel science. The best scientists propounded not only
indispensable theories and addressed not only indispensable problems.
But only indispensable theorems and  problems make these scientists great.

A good theory enables us to settle some indispensable problems. We know many
classical examples of fruitful and powerful theories. Euclidean
geometry and differential calculus were gigantic breakthroughs in the
understanding and mastering the reality. Centuries witness the
strength and power of these theories yielding everyday's solutions of
uncountably many practical problems.
Solution of an indispensable problem is a grind stone for a good
theory since it requires a new conceptual technique and revision of the available
theoretical gadgets. Squaring the circle, the variational principles
of mechanics, and the majority of the Hilbert problems provide
examples of the questions that brought about sweeping changes in
the theoretical outlooks of science.

We must not narrow and simplify the concept of a problem. Science
endeavors to make the complex the simple. Therefore, always actual
are the reconsideration and inventory of the available theories as
well as their simplification, generalization, and unification. The
history of science knows many examples of the perfection, beauty, and
practical power of the theories that arose by way of abstraction and
codification of the preceding views. The success of a new theory
proves that this theory was indispensable.

Freedom in science is the consciousness and appreciation of  the indispensable,
a~vaccine against banality. The call of freedom is inseparable from
the call of mathematics.

\end{document}